\documentclass[a4paper]{article}
\usepackage{mathbbold}
\usepackage{amsfonts}
\usepackage{amssymb}
\usepackage{mathrsfs}
\usepackage{float}
\usepackage{latexcad}
\usepackage{amsmath}
\usepackage{amsfonts,amssymb}
\usepackage{graphicx,color}
\usepackage{psfig}
\usepackage{mathrsfs}

\newtheorem{Theorem}{Theorem}[section]
\newtheorem{Definition}{Definition}[section]
\newtheorem{Lemma}{Lemma}[section]

\catcode`\@=11 \@addtoreset{equation}{section}
\def\theequation{\thesection.\arabic{equation}}
\catcode`\@=12

\def\2{{I \hskip -1.0mm I}}
\def\3{{I \hskip -1.0mm I\hskip -1.0mm I}}
\def\4{{I \hskip -0.9mm V}}
\def\6{{V \hskip -1.35mm I}}

\date{ }

\title{Hyperbolic Mean Curvature Flow}

\author{Chun-Lei He\footnote{Department of Mathematics, Shanghai Jiao Tong University,
Shanghai 200240, China;}\,\,, $\quad$ $\quad$De-Xing
Kong\footnote{Center of Mathematical Sciences, Zhejiang University,
Hangzhou 310027, China;} $\quad$ and $\quad$ Kefeng Liu\footnote{
Department of Mathematics, University of California at Los Angeles,
CA 90095, USA.}\\}
\date{ }
\begin{document}
\maketitle
\begin{abstract}
In this paper we introduce the hyperbolic mean curvature flow and
prove that the corresponding system of partial differential
equations are strictly hyperbolic, and based on this, we show that
this flow admits a unique short-time smooth solution and possesses
the nonlinear stability defined on the Euclidean space with
dimension larger than 4. We derive nonlinear wave equations
satisfied by some geometric quantities related to the hyperbolic
mean curvature flow. Moreover, we also discuss the relation between
the equations for hyperbolic mean curvature flow and the equations
for extremal surfaces in the Minkowski space-time.

\vskip 6mm


\vskip 3mm

\end{abstract}
\newpage
\baselineskip=7mm

\section{Introduction}

Classical differential geometry is on the study of curved spaces and
shapes, in which the time in general does not play a role. However,
in the last few decades, mathematicians have made great strides in
understanding shapes that evolve in time. There are many processes
by which a curve or surface or manifold can evolve, among them two
successful examples are the mean curvature flow and the Ricci flow.
For the Ricci flow, there are many deep and outstanding works, for
example, it can be used to successfully solve the Poincar$\acute{e}$
conjecture and geometrization conjectures. In this paper we will
focus on the mean curvature flow.

It is well known that the mean curvature flow is related on the
motion of surfaces or manifolds. Much more well-known motion of
surfaces are those equating the velocity $\frac{dX}{dt}$ with some
scalar multiple of the normal of the surface. The scalar can be the
curvature, mean curvature or the inverse of the mean curvature with
suitable sign attached. This is the traditional mean curvature flow.
For the traditional mean curvature flow, a beautiful theory has been
developed by Hamilton, Huisken and other researchers (e.g.,
\cite{gh}, \cite{h}, \cite{hu1}), and some important applications
have been obtained, for example, Huisken and Ilmanen developed a
theory of weak solutions of the inverse mean curvature flow and used
it to prove successfully the Riemannian Penrose inequality (see
\cite{hu2}).

A natural problem is as follows: in the above argument if we replace
the velocity $\frac{dX}{dt}$ by the acceleration
$\frac{d^2X}{dt^2}$, what happens? In fact, Yau in \cite{Y} has
suggested the following equation related to a vibrating membrane or
the motion of a surface
\begin{equation}
\frac{d^2X}{dt^2}=H\vec{n},
\end{equation}
where $H$ is the mean curvature and the $\vec{n}$ is the unit inner
normal vector of the surface, and pointed out that very little about
the global time behavior of the hypersurfaces (see page 242 in
\cite{Y}). In deed, according to the authors' knowledge, up to now
only a few of the results on this aspect have been known: a
hyperbolic theory for the evolution of the plane curves has been
developed by Gurtin and Podio-Guidugli \cite{g}, and some
applications to the crystal interfaces have been obtained (see
\cite{R}).

Here we would like to point out that the traditional mean curvature
flow equation is parabolic, however the equation (1.1) is hyperbolic
(see Section 2 for the details). Therefore, in this sense, we name
the equation (1.1) as the hyperbolic version of mean curvature flow,
or {\it hyperbolic mean curvature flow}. Analogous to our recent
work \cite{kl}, in which we introduced and studied the hyperbolic
version of the Ricci flow --- the hyperbolic geometric flow, in this
paper we will investigate the hyperbolic mean curvature flow.

The paper is organized as follows. In Section 2, we introduce the
hyperbolic mean curvature flow and give the short-time existence
theorem. In Section 3, we construct some exact solutions to the
hyperbolic mean curvature flow, these solutions play an important
role in applied fields. Section 4 is devoted to the study on the
nonlinear stability of the hyperbolic mean curvature flow defined on
the Euclidean space with the dimension larger than 4. In Section 5,
we derive the nonlinear wave equations satisfied by some geometric
quantities of the hypersurface $X(\cdot,t)$, these equations show
the wave character of the curvatures. In Section 6, we illustrate
the relations between the hyperbolic mean curvature flow and the
equations for extremal surfaces in the Minkowski space
$\mathbb{R}^{1,n}$.

\section{Hyperbolic mean curvature flow}
Let $\mathscr{M}$ be an $n$-dimensional smooth manifold and
$$X(\cdot,t): \quad \mathscr{M} \rightarrow \mathbb{R}^{n+1}$$ be a one-parameter
family of smooth hypersuface immersions in $\mathbb{R}^{n+1}.$ We
say that it is a solution of the  hyperbolic mean curvature flow if
\begin{equation}\label{2.1}\dfrac{\partial^2 }{\partial
t^2}X(x,t)=H(x,t) \vec{n}(x,t),\qquad \forall ~x\in\mathscr{M}
,\qquad
 \forall ~t>0,\end{equation}where $H(x,t)$ is the mean curvature of $X(x,t)$ and $\vec n(x,t)$ is
the unit inner normal vector on $X(\cdot, t).$

Let $g=\{g_{ij}\}$ and $A=\{h_{ij}\}$ be the induced metric  and the
second fundamental form  on $\mathscr{M}$ in a local coordinate
system $\{x^i\}~ (1\leqslant i\leqslant n)$, respectively.
 Thus, the mean curvature
$H(x,t)$ reads
$$H=g^{ij}h_{ij}.$$
Recall that the Gauss-Weingarten relations
\begin{equation*}\begin{aligned}\dfrac{\partial^2 X}{\partial x^i\partial x^j}=\Gamma^k_{ij}\dfrac{\partial
X}{\partial x^k}+h_{ij}\vec n,~~~~~~ \dfrac{\partial\vec n}{\partial
x^j}=-h_{jl}g^{lm}\dfrac{\partial X}{\partial
x^m}.\end{aligned}\end{equation*} Thus, we have
$$\bigtriangleup_gX=g^{ij}\nabla_i\nabla_jX=g^{ij}\Big(\dfrac{\partial^2X}{\partial x^i\partial x^j}-
\Gamma^k_{ij}\dfrac{\partial X}{\partial x^k}\Big)=g^{ij}h_{ij}\vec
n=H\vec n.$$So the hyperbolic mean curvature flow equation
(\ref{2.1}) can be equivalently rewritten as
\begin{equation}\label{2.2}\dfrac{\partial^2X}{\partial
t^2}=\bigtriangleup_gX=g^{ij}\Big(\dfrac{\partial^2X}{\partial
x^i\partial x^j}- \Gamma^k_{ij}\dfrac{\partial X}{\partial
x^k}\Big).\end{equation} Noting
$$\Gamma^k_{ij}=g^{kl}\Big(\dfrac{\partial^2X}{\partial x^i\partial x^j},\dfrac{\partial X}{\partial x^l}\Big),$$
we get
\begin{equation}\label{2.3}\dfrac{\partial^2X}{\partial
t^2}=g^{ij}\dfrac{\partial^2X}{\partial x^i\partial x^j}-
g^{ij}g^{kl}\Big(\dfrac{\partial^2X}{\partial x^i\partial
x^j},\dfrac{\partial X}{\partial x^l}\Big)\dfrac{\partial
X}{\partial x^k}.\end{equation} It is easy to see that the equation
(\ref{2.3}) is not strictly hyperbolic. Therefore, instead of
considering the equation (\ref{2.3}) we will follow a trick of
DeTurck \cite{D} by modifying the flow through a diffeomorphism of
$\mathscr{M}$, under which (2.3) turns out to be strictly
hyperbolic, so that we can apply the standard theory of hyperbolic
equations.

Suppose $\hat X(x,t) $ is a solution of the equation (\ref{2.1})
(equivalently, (\ref{2.2})) and $\varphi_t: \mathscr{M}\rightarrow
\mathscr{M}$ is a family of diffeomorphisms of $\mathscr{M}$. Let
$$X(x,t)=\varphi^*_t\hat X(x,t),$$ where $\varphi^*_t$ is the pull-back operator of $\varphi_t.$ We
now want to find the evolution equation for the metric $X(x,t).$

Denote
$$y(x,t)=\varphi_t(x)=\{y^1(x,t), y^2(x,t),\cdots, y^n(x,t)\}$$ in
local coordinates, and define $y(x,t)=\varphi_t(x)$ by the following
initial value problem
\begin{equation}\label{2.4}\begin{cases}\dfrac{\partial^2 y^\alpha}{\partial
t^2}=\dfrac{\partial y^\alpha}{\partial x^k}
(g^{jl}(\Gamma^k_{jl}-\tilde{\Gamma}^k_{jl})),\\
y^\alpha(x,0)=x^\alpha,\quad
y^\alpha_t(x,0)=0,\end{cases}\end{equation}where
$\tilde{\Gamma}^k_{jl}$ is the connection corresponding to the
initial metric $ \tilde g_{ij}(x)$. Since
$$\Gamma^k_{jl}=\dfrac{\partial y^\alpha}{\partial
x^j}\dfrac{\partial y^\beta}{\partial x^l}\dfrac{\partial
x^k}{\partial
y^\gamma}\hat\Gamma^\gamma_{\alpha\beta}+\dfrac{\partial
x^k}{\partial y^\alpha}\dfrac{\partial^2y^\alpha}{\partial
x^j\partial x^l},$$the initial value problem (\ref{2.4}) can be
rewritten as
\begin{equation}\label{2.5}\begin{cases}\dfrac{\partial^2 y^\alpha}{\partial
t^2}=g^{jl}\Big(\dfrac{\partial^2 y^\alpha}{\partial x^j\partial
x^l} +\dfrac{\partial y^\beta}{\partial x^j}\dfrac{\partial
y^\gamma}{\partial
x^l}\hat{\Gamma}^\alpha_{\beta\gamma}-\dfrac{\partial
y^\alpha}{\partial x^k}\tilde {\Gamma}^k_{jl}
\Big),\\
y^\alpha(x,0)=x^\alpha,\quad
y^\alpha_t(x,0)=0.\end{cases}\end{equation} Obviously, (\ref{2.5})
is an initial value problem for a strictly hyperbolic system.

On the other hand, noting
\begin{equation*}\begin{aligned}\bigtriangleup_{\hat g}\hat X&=\hat g^{\alpha\beta}\nabla_\alpha\nabla_\beta\hat X
=\hat g^{\alpha\beta}\Big(\dfrac{\partial^2\hat X}{\partial
y^\alpha\partial y^\beta}-\dfrac{\partial \hat X}{\partial
y^\gamma}\hat\Gamma^\gamma_{\alpha\beta}\Big)\\
&=g^{kl}\dfrac{\partial y^\alpha}{\partial x^k}\dfrac{\partial
y^\beta}{\partial x^l}\Big(\dfrac{\partial}{\partial
y^\alpha}\Big(\dfrac{\partial X}{\partial x^i}\dfrac{\partial
x^i}{\partial y^\beta}\Big)-\dfrac{\partial X}{\partial
x^i}\dfrac{\partial x^i}{\partial
y^\gamma}\hat\Gamma^\gamma_{\alpha\beta}\Big)\\
&=g^{kl}\dfrac{\partial^2X}{\partial x^k\partial
x^l}+g^{kl}\dfrac{\partial y^\alpha}{\partial x^k}\dfrac{\partial
y^\beta}{\partial x^l}\dfrac{\partial X}{\partial
x^i}\dfrac{\partial^2x^i}{\partial y^\alpha\partial
y^\beta}-g^{kl}\dfrac{\partial X}{\partial
x^i}\Big(\Gamma^i_{kl}-\dfrac{\partial x^i}{\partial
y^\gamma}\dfrac{\partial^2y^\gamma}{\partial x^k\partial x^l}\Big)\\
&=g^{kl}\nabla_k\nabla_lX=\bigtriangleup
_gX,\end{aligned}\end{equation*} we have
\begin{equation*}\begin{aligned}\dfrac{\partial^2X}{\partial t^2}&=\dfrac{\partial^2\hat X}
{\partial y^\alpha\partial y^\beta}\dfrac{\partial
y^\alpha}{\partial t}\dfrac{\partial y^\beta}{\partial
t}+2\dfrac{\partial^2\hat X} {\partial t\partial
y^\beta}\dfrac{\partial y^\beta}{\partial t}+\dfrac{\partial^2\hat
X} {\partial t^2}+\dfrac{\partial\hat X} {\partial
y^\alpha}\dfrac{\partial^2 y^\alpha}{\partial
t^2}\\
&=\dfrac{\partial^2\hat X} {\partial y^\alpha\partial
y^\beta}\dfrac{\partial y^\alpha}{\partial t}\dfrac{\partial
y^\beta}{\partial t}+2\dfrac{\partial^2\hat X} {\partial t\partial
y^\beta}\dfrac{\partial y^\beta}{\partial t}+\bigtriangleup_{\hat
g}\hat X+\dfrac{\partial X}{\partial x^k}\dfrac{\partial x^k}
{\partial y^\alpha}\dfrac{\partial^2 y^\alpha}{\partial
t^2}\\
&=\bigtriangleup_gX+\dfrac{\partial X}{\partial
x^k}g^{ij}(\Gamma^k_{ij}-\tilde{\Gamma}^k_{ij})+\dfrac{\partial^2\hat
X} {\partial y^\alpha\partial y^\beta}\dfrac{\partial
y^\alpha}{\partial t}\dfrac{\partial y^\beta}{\partial
t}+2\dfrac{\partial^2\hat X} {\partial t\partial
y^\beta}\dfrac{\partial y^\beta}{\partial t}\\
&=g^{ij}\Big(\dfrac{\partial^2 X}{\partial x^i\partial
x^j}-\dfrac{\partial X}{\partial
x^k}\Gamma^k_{ij}\Big)+\dfrac{\partial X}{\partial
x^k}g^{ij}(\Gamma^k_{ij}-\tilde{\Gamma}^k_{ij})+\dfrac{\partial^2\hat
X} {\partial y^\alpha\partial y^\beta}\dfrac{\partial
y^\alpha}{\partial t}\dfrac{\partial y^\beta}{\partial
t}+2\dfrac{\partial^2\hat X} {\partial t\partial
y^\beta}\dfrac{\partial y^\beta}{\partial t}\\
&=g^{ij}\dfrac{\partial^2 X}{\partial x^i\partial
x^j}-g^{ij}\tilde{\Gamma}^k_{ij}\dfrac{\partial X}{\partial
x^k}+\dfrac{\partial^2\hat X} {\partial y^\alpha\partial
y^\beta}\dfrac{\partial y^\alpha}{\partial t}\dfrac{\partial
y^\beta}{\partial t}+2\dfrac{\partial^2\hat X} {\partial t\partial
y^\beta}\dfrac{\partial y^\beta}{\partial
t}.\end{aligned}\end{equation*} By the standard theory of hyperbolic
equations (see \cite{H}), we have the following result.

\begin{Theorem}\label{t1} {\bf (Local existences and uniqueness)} Let $\mathscr{M}$
 be an n-dimensional smooth compact manifold, and $X_0$ be a
smooth hypersuface immersion of $\mathscr{M}$ into
$\mathbb{R}^{n+1}.$ Then there exists a constant $T
> 0$ such that the initial value problem
\begin{equation}\label{2.6}\begin{cases} \dfrac{\partial^2 }{\partial
t^2}X(x,t)=H(x,t)\vec n(x,t),\\
X\big |_{t=0}=X_0(x), ~~\dfrac{\partial X}{\partial t}(x,t)\Big
|_{t=0}=X_1(x)\vspace{2mm}\end{cases}\end{equation} has a unique
smooth solution $X(x,t)$ on $\mathscr{M}\times[0,T),$ where $X_1(x)$
is a smooth vector-valued function on $\mathscr{M}.$
\end{Theorem}

\section{Exact solutions}
In order to understand further the hyperbolic mean curvature flow,
in this section we investigate some exact solutions. These exact
solutions play an important role in applied fields. To do so, we
first consider the following initial value problem for an  ordinary
differential equation
\begin{equation}\label{3.1}\begin{cases}r_{tt}=-\dfrac{1}{r}~,\\
r(0)=r_0>0~,~~~r_t(0)=r_1~.\end{cases}\end{equation} For this
initial value problem, we have the following lemma.
\begin{Lemma} For arbitrary initial data $r_0>0$, if the initial velocity  $r_1\leqslant
0$, then the solution $r=r(t)$  decreases and attains its zero point
at time $t_0$ (in particular, when $r_1=0$, we have
$t_0=\sqrt{\frac{\pi}{2}}r_0$); if the initial velocity is positive,
then the solution $r$ increases first and then decreases and attains
its zero point in a finite time $t_0$.
\end{Lemma}
{\bf Proof.} ~The proof is similar to the arguments in \cite{R}. The
 following discussion is divided into two cases.

 {\bf Case \textrm{I}} $\quad$ The initial velocity is
nonpositive, i.e., $r_1\leqslant 0.$

We argue by contradiction. Let us assume that $r(t)>0$ for all time
$t>0.$ Then $r_{tt}<0$ and $r_t(t)<r_t(0)=r_1\leqslant 0$ for $t>0.$
Hence there exists a time $t_0$ such that $r(t_0)=0$ (see Fig.1).
This is a contradiction.

Moreover, for the case $r_1=0,$ we can derive the explicit
expression for $t_0$ according to the equation (\ref{3.1}).
Multiplying both side of $r_{tt}=-\dfrac1r$ by $r_t,$ integrating,
applying the initial condition $r_t(0)\leqslant 0$ and $r(t_0)=0$,
integrating once again yields$$\int_0^{t_0}\dfrac{\sqrt
2}{2r_0}dt=\int_0^{\infty}e^{-u^2}du=\dfrac{\sqrt \pi}{2},$$ where
$u=\sqrt{\ln\dfrac{r_0}{r}}.$ Thus we obtain
$$t_0=\sqrt{\dfrac{\pi}{2}}r_0.$$ \vskip 2mm

\begin{figure}[H]
    \begin{center}
\begin{picture}(328,164)
\thinlines \drawvector{4.0}{40.0}{146.0}{1}{0}
\drawvector{22.0}{22.0}{138.0}{0}{1}
\drawvector{172.0}{40.0}{152.0}{1}{0}
\drawvector{194.0}{22.0}{138.0}{0}{1}
\path(22.0,112.0)(22.0,112.0)(23.43,111.98)(24.84,111.95)(26.27,111.91)(27.66,111.86)(29.06,111.79)(30.43,111.69)(31.81,111.58)(33.16,111.47)
\path(33.16,111.47)(34.52,111.33)(35.86,111.18)(37.18,111.0)(38.5,110.8)(39.79,110.61)(41.09,110.38)(42.38,110.15)(43.65,109.9)(44.9,109.62)
\path(44.9,109.62)(46.16,109.33)(47.4,109.02)(48.63,108.72)(49.84,108.37)(51.06,108.02)(52.26,107.66)(53.44,107.26)(54.62,106.87)(55.77,106.44)
\path(55.77,106.44)(56.94,106.01)(58.08,105.56)(59.2,105.09)(60.33,104.62)(61.44,104.11)(62.55,103.59)(63.62,103.06)(64.7,102.51)(65.77,101.94)
\path(65.77,101.94)(66.83,101.37)(67.87,100.76)(68.91,100.15)(69.94,99.51)(70.94,98.87)(71.95,98.2)(72.94,97.52)(73.93,96.83)(74.9,96.12)
\path(74.9,96.12)(75.86,95.38)(76.8,94.63)(77.75,93.87)(78.66,93.09)(79.58,92.3)(80.48,91.5)(81.38,90.66)(82.26,89.81)(83.15,88.95)
\path(83.15,88.95)(84.01,88.08)(84.86,87.19)(85.69,86.27)(86.52,85.34)(87.34,84.41)(88.16,83.44)(88.94,82.48)(89.73,81.48)(90.51,80.47)
\path(90.51,80.47)(91.27,79.44)(92.04,78.41)(92.77,77.34)(93.51,76.27)(94.23,75.19)(94.94,74.08)(95.65,72.94)(96.33,71.81)(97.01,70.66)
\path(97.01,70.66)(97.68,69.48)(98.33,68.3)(98.98,67.08)(99.62,65.87)(100.23,64.62)(100.86,63.37)(101.45,62.11)(102.05,60.81)(102.62,59.52)
\path(102.62,59.52)(103.2,58.18)(103.76,56.86)(104.3,55.5)(104.84,54.13)(105.37,52.75)(105.9,51.34)(106.4,49.93)(106.9,48.49)(107.37,47.04)
\path(107.37,47.04)(107.84,45.58)(108.31,44.09)(108.76,42.59)(109.2,41.06)(109.63,39.54)(110.05,37.99)(110.47,36.41)(110.87,34.84)(111.25,33.24)
\path(111.25,33.24)(111.62,31.63)(111.98,30.0)(112.0,30.0)
\path(194.0,88.0)(194.0,88.0)(195.14,90.25)(196.3,92.44)(197.47,94.58)(198.63,96.66)(199.77,98.68)(200.92,100.65)(202.08,102.55)(203.24,104.4)
\path(203.24,104.4)(204.38,106.19)(205.53,107.94)(206.67,109.61)(207.83,111.23)(208.97,112.8)(210.11,114.3)(211.25,115.76)(212.39,117.15)(213.53,118.48)
\path(213.53,118.48)(214.67,119.76)(215.82,120.98)(216.96,122.16)(218.08,123.26)(219.22,124.3)(220.36,125.3)(221.49,126.24)(222.61,127.11)(223.75,127.94)
\path(223.75,127.94)(224.88,128.71)(226.0,129.41)(227.13,130.05)(228.25,130.66)(229.38,131.19)(230.5,131.66)(231.61,132.08)(232.74,132.44)(233.86,132.75)
\path(233.86,132.75)(234.97,133.0)(236.08,133.19)(237.21,133.33)(238.32,133.41)(239.42,133.44)(240.55,133.39)(241.66,133.3)(242.77,133.14)(243.86,132.94)
\path(243.86,132.94)(244.98,132.67)(246.08,132.36)(247.19,131.97)(248.29,131.53)(249.39,131.05)(250.48,130.5)(251.58,129.88)(252.69,129.22)(253.79,128.5)
\path(253.79,128.5)(254.89,127.72)(255.98,126.88)(257.07,125.98)(258.17,125.02)(259.26,124.01)(260.35,122.95)(261.42,121.83)(262.51,120.65)(263.6,119.41)
\path(263.6,119.41)(264.69,118.12)(265.77,116.76)(266.85,115.36)(267.94,113.88)(269.01,112.37)(270.1,110.79)(271.17,109.15)(272.25,107.45)(273.32,105.69)
\path(273.32,105.69)(274.39,103.88)(275.48,102.02)(276.54,100.09)(277.61,98.12)(278.69,96.08)(279.76,93.98)(280.82,91.83)(281.89,89.62)(282.95,87.36)
\path(282.95,87.36)(284.01,85.02)(285.07,82.65)(286.14,80.2)(287.2,77.7)(288.26,75.16)(289.32,72.55)(290.36,69.87)(291.42,67.16)(292.48,64.37)
\path(292.48,64.37)(293.52,61.54)(294.58,58.63)(295.64,55.68)(296.69,52.66)(297.73,49.61)(298.77,46.47)(299.82,43.29)(300.86,40.06)(301.91,36.75)
\path(301.91,36.75)(302.95,33.4)(303.98,30.0)(304.0,30.0)
\drawcenteredtext{16.5}{158.0}{$r$}
\drawcenteredtext{148.0}{32.0}{$t$}
\drawcenteredtext{322.0}{32.0}{$t$}
\drawcenteredtext{188.0}{158.0}{$r$}
\drawcenteredtext{16.0}{32.0}{$0$}
\drawcenteredtext{105.0}{32.0}{$t_0$}
\drawcenteredtext{296.0}{32.0}{$t_0$}
\drawcenteredtext{186.0}{32.0}{$0$}
\drawcenteredtext{16.0}{112.0}{$r_0$}
\drawcenteredtext{188.0}{88.0}{$r_0$}
\drawcenteredtext{92.0}{102.0}{$r=r(t)$}
\drawcenteredtext{286.0}{122.0}{$r=r(t)$}
\drawcenteredtext{68.0}{0.0}{Fig. 1: $r_1\leqslant0$}
\drawcenteredtext{240.0}{0.0}{Fig. 2: $r_1>0$}
\end{picture}
    \end{center}
\end{figure}
 {\bf Case \textrm{II}}$\quad$ The initial velocity is positive,
i.e., $r_1>0.$

By (\ref{3.1}), we obtain
$$r_t^2=-2\ln r+2\ln r_0+r_1^2~.$$ Then we have
$$r\leqslant e^{\frac{r_1^2}{2}}r_0~.$$ If $r$ increases for all time, i.e., $r_t>0$ for all time $t$, we have
$r_0<r\leqslant e^{\frac{r_1^2}{2}}r_0$ and
$-\frac1{r_0}<r_{tt}\leqslant -e^{-\frac{r_1^2}{2}}\frac1{r_0}.$
Thus, the curve $r_t$ can be bounded by two straight lines
$r_t=-\dfrac{1}{r_0}t+r_1$ and
$r_t=-\frac1{r_0}e^{-\frac{r_1^2}{2}}t+r_1$. On the other hand,
$r_t$ is a convex function since $$(r_t)_{tt}=\dfrac{r_t}{r^2}>0.$$
Therefore, $r_t$ will change sign and becomes negative at certain
finite time, this contradicts to the hypothesis that $r$ is always
increasing. Thus, in this case, $r$ increases first and then
decreases and attains its zero point in a finite time (see Fig. 2).
The proof is finished. $\quad\quad\blacksquare$

\vskip 5mm

In what follows, we are interested in some exact solutions of
hyperbolic mean curvature flow (\ref{2.1}).

{\bf Example 1:} Consider a family of  spheres
\begin{equation}\label{3.2}X(x,t)=r(t)(\cos\alpha\cos\beta,\cos\alpha\sin\beta,\sin\alpha),\end{equation}where
$\alpha\in [-\dfrac{\pi}{2},\dfrac{\pi}{2}],~ \beta\in [0,2\pi]$.

 Clearly, the induced metric and the second fundamental form are,
 respectively,
$$g_{11}=r^2,~~ g_{22}=r^2\cos^2\alpha,~~ g_{12}=g_{21}=0$$and
$$h_{11}=r,~~ h_{22}=r\cos^2\alpha,~~ h_{12}=h_{21}=0.$$ The mean
curvature is $$H=\dfrac{2}{r}.$$ On the other hand, the Christoffel
symbols read
\begin{equation*}\begin{aligned}&\Gamma_{11}^{1}=\Gamma_{12}^{1}=0,\qquad\Gamma^{1}_{22}=\cos\alpha\sin\alpha,\\
&\Gamma_{11}^{2}=\Gamma^{2}_{22}=0,\qquad\Gamma^{2}_{12}=-\dfrac{\sin\alpha}{\cos\alpha}.
\end{aligned}\end{equation*}
Thus, we obtain from (\ref{2.1}) or (\ref{2.2}) that
\begin{equation}\label{32}r_{tt}=-\dfrac2r.\end{equation}

By Lemma 3.1, it can be easily observed that, for arbitrary
$r(0)>0,$ if $ r_t(0)\leqslant0$, the evolving sphere will shrink to
a point; if $r_t(0)>0,$ the evolving sphere will expand first and
then shrink to a point.$\quad\quad\blacksquare$

In fact, this phenomena can also be interpreted by physical
principles. From (\ref{3.2}), we have
\begin{equation}\label{3.4}X_t(x,t)=r_t(t)(\cos\alpha\cos\beta,\cos\alpha\sin\beta,\sin\alpha)\end{equation}and
\begin{equation}\label{3.5}X_{tt}(x,t)=r_{tt}(t)(\cos\alpha\cos\beta,\cos\alpha\sin\beta,\sin\alpha).\end{equation}
By (\ref{32}) and (\ref{3.5}), the direction of acceleration is
always the same as the inner normal vector. Thus, due to
(\ref{3.4}), if $ r_t(0)\leqslant0$, i.e., the initial velocity
direction is the same as inner normal vector, then evolving sphere
will shrink to a point; if $r_t(0)>0$, i.e., the initial velocity
direction is opposite to inner normal vector, then the evolving
sphere will expand first and then shrink to a point.

{\bf Example 2:}  We now consider an exact solution with axial
symmetry. In other words, we focus on the cylinder solution for the
hyperbolic mean curvature flow  which takes the following form
\begin{equation*}X(x,t)=(r(t)\cos\alpha,~r(t)\sin\alpha,~\rho),\end{equation*}where
$\alpha\in [0,2\pi],~\rho\in[0,\rho_0]$.

Obviously, the induced metric and the second fundamental form read,
 respectively,
$$g_{11}=r^2,~~ g_{22}=1,~~ g_{12}=g_{21}=0$$and$$h_{11}=r,~~ h_{22}=0,~~
h_{12}=h_{21}=0.$$ The mean curvature is $$H=\dfrac1r.$$ Moreover,
the Christoffel symbols are
\begin{equation*}\Gamma_{ij}^{k}=0,\quad\quad \forall~ i, j, k=1,2.
\end{equation*} Then, we obtain from (\ref{2.1}) or (\ref{2.2}) that
$$r_{tt}=-\dfrac1r.$$

By Lemma 3.1, it can be easily found that the evolving cylinder will
always shrink to a straight line for arbitrary $\rho_0>0, r(0)>0$
and $ r_t(0).$$\quad\quad\blacksquare$

\section{Nonlinear stability}

In this section, we consider the nonlinear stability of the
hyperbolic mean curvature flow defined on the Euclidean space with
the dimension larger than 4.

Let $\mathscr{M}$ be an $n$-dimensional ($n>4$) complete Riemannian
manifold. Given the hypersurfaces $X_1(x)$ and $X_2(x)$ on $\mathscr
M,$ we consider the following initial value problem
\begin{equation}\label{4.1}\begin{cases} \dfrac{\partial^2 }{\partial
t^2}X(x,t)=H(x,t)\vec n(x,t),\\
X(x,0)=X_0(x)+\varepsilon X_1(x),~~~\dfrac{\partial X}{\partial
t}(x,0)=\varepsilon
X_2(x)\vspace{2mm},\end{cases}\end{equation}where $\varepsilon>0$ is
a small parameter.

A coordinate chart $(x^1,\cdots,x^n)$ on a Riemannian manifold
$(\mathscr{M},g)$ is called {\it harmonic} if $$\bigtriangleup x^j=0
\quad (j=1,\cdots,n).$$ DeTurck \cite{D} showed that a coordinate
function $x^k$ is {\it harmonic} if and only if
$$\bigtriangleup
x^k=g^{ij}\dfrac{\partial^2 x^k}{\partial x^i\partial
x^j}-g^{ij}\Gamma_{ij}^{l}\dfrac{\partial x^k}{\partial
x^l}=-\Gamma^k_{ij}g^{ij}=-\Gamma^k=0.$$  He also proved the
following theorem on the existence of harmonic coordinates.
\begin{Lemma}\label{l1} Let the metric $g$ on a Riemannian manifold $(\mathscr{M},g)$ be of
class $C^{k,\alpha}$ (for $k\geqslant1)$ (resp. $C^{\omega}$) in a
local coordinate chart about some point $p.$ Then there is a
neighborhood of $p$ in which harmonic coordinates exist, these new
coordinates being $C^{k+1,\alpha}$ (resp. $C^{\omega}$) functions of
the original coordinates. Moreover, all harmonic coordinate charts
defined near $p$ have this regularity.\end{Lemma}

By Lemma \ref{l1} and Theorem \ref{t1}, we can choose the harmonic
coordinates around a fixed point $p\in \mathscr{M}$ and for a fixed
time $t\in \mathbb{R}^+.$ Then the hyperbolic mean curvature flow
(\ref{2.2}) can be equivalent by  written as
\begin{equation*}\dfrac{\partial^2X}{\partial t^2}=g^{ij}\dfrac{\partial^2X}{\partial x^i\partial x^j}.
\end{equation*}

\begin{Definition}$X_0(x)$ possesses the (locally) nonlinear stability with respect to
$(X_1(x),X_2(x)),$ if there exists a positive constant
$\varepsilon_0=\varepsilon_0(X_1(x),X_2(x))$ such that, for any
$\varepsilon\in (0,\varepsilon_0],$ the initial value problem
(\ref{4.1}) has a unique (local) smooth solution $X(x,t);$

$X_0(x)$ is said to be (locally) nonlinear stable if it possesses
the (locally) nonlinear stability with respect to arbitrary $X_1(x)$
and $X_2(x).$
\end{Definition}

\begin{Theorem}$X_0(x)=(x^1,x^2,\cdots,x^n,0)$ ( $n>4$) is nonlinearly stable.\end{Theorem}
{\bf Proof.} Choose the harmonic coordinates around a fixed point
$p\in \mathscr{M}$ and for a fixed time $t\in \mathbb{R}^+.$ Then
the equation (\ref{4.1}) can be written as
\begin{equation}\label{2.7}\begin{cases} \dfrac{\partial^2 }{\partial
t^2}X(x,t)=g^{ij}\dfrac{\partial^2X}{\partial x^i\partial x^j},\\
X(x,0)=X_0(x)+\varepsilon X_1(x),~~~\dfrac{\partial X}{\partial
t}(x,0)=\varepsilon X_2(x).\end{cases}\end{equation}

Define $Y(x,t)=(y^1,\cdots,y^n,y^{n+1})$ in the following way
$$X(x,t)=X_0(x)+Y(x,t).$$ Then for small $Y(x,t),$ we have
\begin{equation}\label{2.8}g_{ij}=\Big(\dfrac{\partial X}{\partial x^i},
\dfrac{\partial X}{\partial x^j}\Big)=\delta_{ij}+\dfrac{\partial
y^j}{\partial x^i}+\dfrac{\partial y^i}{\partial
x^j}+y_{ij}~\end{equation} and
\begin{equation}\label{2.9}g^{ij}=\delta^{ij}-\dfrac{\partial
y^j}{\partial x^i}-\dfrac{\partial y^i}{\partial
x^j}-y_{ij}+O(||\lambda||^2),\end{equation} where
$$\delta_{ij}=\begin{cases}1,~~~i=j~(i,j=1,\cdots,n),\\
0,~~~i\neq j~(i,j=1,\cdots,n),\end{cases}$$
$$y_{ij}=\Big(\dfrac{\partial Y}{\partial x^i}, \dfrac{\partial
Y}{\partial x^j}\Big)=~\sum_{p=1}^{n+1}\dfrac{\partial y^p}{\partial
x^i}\dfrac{\partial y^p}{\partial x^j},$$
$$\lambda=\left(\dfrac{\partial y^p}{\partial
x^q}\right)\quad (p=1,2,\cdots,n+1;~~q=1,2,\cdots,n).$$
 The equation
(\ref{2.7}) can be rewritten as
\begin{equation}\label{2.10}\begin{cases} \dfrac{\partial^2 }{\partial
t^2}Y(x,t)=\Big(\delta^{ij}-\dfrac{\partial y^j}{\partial
x^i}-\dfrac{\partial y^i}{\partial
x^j}-y_{ij}+O(||\lambda||^2)\Big)\dfrac{\partial^2Y}{\partial
x^i\partial x^j},\vspace{2mm}\\ Y(x,0)=\varepsilon
X_1(x),~~~\dfrac{\partial Y}{\partial t}(x,0)=\varepsilon
X_2(x).\end{cases}\end{equation}

Define $$\hat\lambda=\Big(\dfrac{\partial y^m}{\partial
x^k},\dfrac{\partial^2 y^m}{\partial x^k\partial
x^l}\Big)~(m=1,\cdots,n+1;~k,l=1,\cdots,n),$$ then for all $p$ we
have
\begin{equation*}\begin{aligned}\dfrac{\partial^2y^p}{\partial t^2}&=\dfrac{\partial^2y^p}{\partial x^i\partial
x^i}+ \Big(-\dfrac{\partial y^j}{\partial x^i}-\dfrac{\partial
y^i}{\partial
x^j}-y_{ij}+O(||\lambda||^2)\Big)\dfrac{\partial^2y^p}{\partial
x^i\partial x^j}\\
&=\dfrac{\partial^2y^p}{\partial x^i\partial
x^i}+O(||\hat\lambda||^2).\end{aligned}\end{equation*} By the
well-known result on the global existence for  nonlinear wave
equations (e.g., see \cite{C}, \cite{H}, \cite{K}), there exists a
unique global smooth solution $Y=Y(x,t)$ for the Cauchy problem
(\ref{2.10}). Thus, the proof is completed.$\quad\quad\blacksquare$

\section{Evolution of metric and curvatures}
From the evolution equation (\ref{2.1}) for the hyperbolic mean
curvature flow, we can derive the evolution equations for some
geometric quantities of the hypersurface $X(\cdot,t)$, these
equations will play an important role in the future study on the
hyperbolic mean curvature flow.

\begin{Lemma} Under the hyperbolic mean curvature flow, the following identities hold \begin{equation}\label{5.1}
\bigtriangleup h_{ij}=\nabla_i\nabla_jH+Hh_{il}g^{lm}h_{mj}-|A|^2h_{ij},
\end{equation}\begin{equation}\label{5.2}\bigtriangleup |A|^2=2g^{ik}g^{jl}h_{kl}\nabla_i\nabla_jH+2|\nabla A|^2
+2H{\bf{tr}}(A^3)-2|A|^4,\end{equation} where
$$|A|^2=g^{ij}g^{kl}h_{ik}h_{jl}~,~~{\bf{tr}}(A^3)=g^{ij}g^{kl}g^{mn}h_{ik}h_{lm}h_{nj}.$$
\end{Lemma}

Lemma 5.1 can be found in Zhu $\cite{Z}$.

\begin{Theorem}Under the hyperbolic mean curvature flow, it holds
that
\begin{equation}\label{5.3}\dfrac{\partial^2 g_{ij}}{\partial t^2}=-2Hh_{ij}+2
\Big(\dfrac{\partial^2X}{\partial t\partial
x^i},\dfrac{\partial^2X}{\partial t\partial x^j}\Big),\end{equation}
\begin{equation}\label{5.4}\dfrac{\partial^2\vec n}{\partial t^2}=-g^{ij}\dfrac{\partial
H}{\partial x^i}\dfrac{\partial X}{\partial x^j} +g^{ij}\Big(\vec
n,\dfrac{\partial^2X}{\partial t\partial
x^i}\Big)\Big[2g^{kl}\Big(\dfrac{\partial X}{\partial
x^j},\dfrac{\partial^2X}{\partial t\partial x^l}\Big)\dfrac{\partial
X}{\partial x^k}+g^{kl}\Big(\dfrac{\partial X}{\partial
x^l},\dfrac{\partial^2X}{\partial t\partial x^j}\Big)\dfrac{\partial
X}{\partial x^k}-\dfrac{\partial^2X}{\partial t\partial x^j}\Big]
\end{equation}
and \begin{equation}\label{5.5}\dfrac{\partial^2h_{ij}}{\partial
t^2}=\bigtriangleup
h_{ij}-2Hh_{il}h_{mj}g^{lm}+|A|^2h_{ij}+g^{kl}h_{ij}\Big(\vec
n,\dfrac{\partial^2X}{\partial t\partial x^k}\Big)\Big(\vec
n,\dfrac{\partial^2X}{\partial t\partial
x^l}\Big)-2\dfrac{\partial\Gamma^k_{ij}}{\partial t}\Big(\vec
n,\dfrac{\partial^2X}{\partial t\partial x^k}\Big).~~~~~
\end{equation}
\end{Theorem}

\noindent {\bf Proof.} With the aids of the definitions of the
metric, the second fundamental form and the Gauss-Weingarten
relations, we can give a complete proof of Theorem 5.1. In fact, by
the definition of the metric, we have
\begin{equation*}\begin{aligned}\dfrac{\partial^2g_{ij}}
{\partial t^2}&=\dfrac{\partial^2}{\partial t^2}
\Big(\dfrac{\partial X}{\partial x^i},\dfrac{\partial X}{\partial
x^j}\Big)=\Big(\dfrac{\partial^3X}{\partial t^2\partial
x^i},\dfrac{\partial X}{\partial
x^j}\Big)+2\Big(\dfrac{\partial^2X}{\partial t\partial
x^i},\dfrac{\partial^2X}{\partial t\partial
x^j}\Big)+(\dfrac{\partial X}{\partial
x^i},\dfrac{\partial^3X}{\partial t^2\partial
x^j})\vspace{5mm}\\
&=\Big(\dfrac{\partial}{\partial x^i}(H\vec n),\dfrac{\partial
X}{\partial x^j}\Big)+2\Big(\dfrac{\partial^2X}{\partial t\partial
x^i},\dfrac{\partial^2X}{\partial t\partial
x^j}\Big)+\Big(\dfrac{\partial X}{\partial
x^i},\dfrac{\partial}{\partial x^j}(H\vec
n)\Big)\vspace{5mm}\\
&=H\Big(-h_{ik}g^{kl}\dfrac{\partial X}{\partial
x^l},\dfrac{\partial X}{\partial
x^j}\Big)+2\Big(\dfrac{\partial^2X}{\partial t\partial
x^i},\dfrac{\partial^2X}{\partial t\partial
x^j}\Big)+H\Big(\dfrac{\partial X}{\partial
x^i},-h_{jk}g^{kl}\dfrac{\partial X}{\partial
x^l}\Big)\vspace{5mm}\\
&=-2Hh_{ij}+2\Big(\dfrac{\partial^2X}{\partial t\partial
x^i},\dfrac{\partial^2X}{\partial t\partial
x^j}\Big).\\\end{aligned}\end{equation*} This gives the proof of
(\ref{5.3}).

On the other hand,
$$\dfrac{\partial\vec n}{\partial t}=\left(\dfrac{\partial\vec
n}{\partial t},\dfrac{\partial X}{\partial
x^i}\right)g^{ij}\dfrac{\partial X}{\partial x^j}=-\left(\vec n,
\dfrac{\partial^2 X}{\partial t\partial
x^i}\right)g^{ij}\dfrac{\partial X}{\partial x^j},$$ then
\begin{equation*}\begin{aligned}\dfrac{\partial^2\vec n}{\partial t^2}
&=-\Big(\dfrac{\partial\vec n}{\partial
t},\dfrac{\partial^2X}{\partial t\partial
x^i}\Big)g^{ij}\dfrac{\partial X}{\partial x^j}-\Big(\vec
n,\dfrac{\partial^3X}{\partial t^2\partial
x^i}\Big)g^{ij}\dfrac{\partial X}{\partial x^j}+\Big(\vec
n,\dfrac{\partial^2X}{\partial t\partial
x^i}\Big)g^{ik}g^{jl}\dfrac{\partial g_{kl}}{\partial
t}\dfrac{\partial X}{\partial x^j}-g^{ij}\Big(\vec
n,\dfrac{\partial^2X}{\partial t\partial x^i}\Big)\dfrac{\partial
^2X}{\partial t\partial
x^j}\\
&=g^{ij}g^{kl}\Big(\vec n,\dfrac{\partial^2X}{\partial t\partial
x^k}\Big)\Big(\dfrac{\partial X}{\partial
x^l},\dfrac{\partial^2X}{\partial t\partial x^i}\Big)\dfrac{\partial
X}{\partial x^j}-\Big(\vec n,\dfrac{\partial }{\partial x^i}(H\vec
n)\Big)g^{ij}\dfrac{\partial
X}{\partial x^j}+\\
&~~~~g^{ik}g^{jl}\Big(\vec n,\dfrac{\partial^2X}{\partial t\partial
x^i}\Big)\Big[\Big(\dfrac{\partial^2X}{\partial t\partial
x^k},\dfrac{\partial X}{\partial x^l}\Big)+\Big(\dfrac{\partial
X}{\partial x^k},\dfrac{\partial^2X}{\partial t\partial
x^l}\Big)\Big]\dfrac{\partial X}{\partial x^j}-\Big(\vec
n,\dfrac{\partial^2X }{\partial t\partial
x^i}\Big)g^{ij}\dfrac{\partial^2 X}{\partial t\partial
x^j}\\
&=-g^{ij}\dfrac{\partial H}{\partial x^i}\dfrac{\partial X}{\partial
x^j}-g^{ij}\Big(\vec n,\dfrac{\partial^2X}{\partial t\partial
x^i}\Big)\dfrac{\partial^2X}{\partial t\partial
x^j}+\\
&~~~~g^{ij}g^{kl}\Big(\vec n,\dfrac{\partial^2X}{\partial t\partial
x^i}\Big)\Big[\Big(\dfrac{\partial^2X}{\partial t\partial
x^j},\dfrac{\partial X}{\partial x^l}\Big)+2\Big(\dfrac{\partial
X}{\partial x^j},\dfrac{\partial^2X}{\partial t\partial
x^l}\Big)\Big]\dfrac{\partial X}{\partial x^k}.
\end{aligned}\end{equation*}This proves (\ref{5.4}).

By virtue of
$$\dfrac{\partial h_{ij}}{\partial t}=\dfrac{\partial}{\partial t}\Big(\vec n,\dfrac{\partial^2X}
{\partial x^i\partial x^j}\Big) =\Big(\dfrac{\partial\vec
n}{\partial t},\dfrac{\partial^2X}{\partial x^i\partial
x^j}\Big)+\Big(\vec n,\dfrac{\partial^3X}{\partial t\partial
x^i\partial x^j}\Big),$$ we have
\begin{equation*}\begin{aligned}\dfrac{\partial^2h_{ij}}{\partial t^2}&=\Big(\dfrac{\partial^2\vec n}{\partial
t^2},\dfrac{\partial^2X}{\partial x^i\partial
x^j}\Big)+2\Big(\dfrac{\partial\vec n}{\partial
t},\dfrac{\partial^3X}{\partial t\partial x^i\partial
x^j}\Big)+\Big(\vec n,\dfrac{\partial^4X}{\partial t^2\partial
x^i\partial
x^j}\Big)\\
&=-g^{kl}\Big(\dfrac{\partial H}{\partial x^k}\dfrac{\partial
X}{\partial x^l},\dfrac{\partial^2X}{\partial x^i\partial
x^j}\Big)-g^{kl}\Big(\vec n,\dfrac{\partial^2X}{\partial t\partial
x^k}\Big)\Big(\dfrac{\partial^2 X}{\partial t\partial
x^l},\dfrac{\partial^2X}{\partial x^i\partial x^j}\Big)\\
&~~+g^{pq}g^{kl}\Big(\vec n,\dfrac{\partial^2X}{\partial t\partial
x^p}\Big)\Big[\Big(\dfrac{\partial X}{\partial
x^l},\dfrac{\partial^2X}{\partial t\partial
x^q}\Big)+2\Big(\dfrac{\partial X}{\partial
x^q},\dfrac{\partial^2X}{\partial t\partial
x^l}\Big)\Big](\dfrac{\partial X}{\partial
x^k},\dfrac{\partial^2X}{\partial x^i\partial x^j})\\
&~~-2g^{kl}\Big(\vec n,\dfrac{\partial^2X}{\partial t\partial
x^k}\Big)\Big(\dfrac{\partial X}{\partial
x^l},\dfrac{\partial^3X}{\partial t\partial x^i\partial
x^j}\Big)+\Big(\vec n,\dfrac{\partial^2}{\partial
x^i\partial x^j}(H\vec n)\Big)\\
&=-\dfrac{\partial H}{\partial
x^k}\Gamma^k_{ij}-g^{kl}\Gamma^m_{ij}\Big(\vec n
,\dfrac{\partial^2X}{\partial t\partial
x^k}\Big)\Big(\dfrac{\partial X}{\partial
x^m},\dfrac{\partial^2X}{\partial t\partial
x^l}\Big)-g^{kl}h_{ij}\Big(\vec n,\dfrac{\partial^2 X}{\partial
t\partial x^k}\Big)\Big(\vec n,\dfrac{\partial^2 X}{\partial
t\partial x^l}\Big)\\
&~~+g^{pq}\Gamma^l_{ij}\Big(\vec n,\dfrac{\partial^2 X}{\partial
t\partial x^p}\Big)\Big[\Big(\dfrac{\partial X}{\partial
x^l},\dfrac{\partial^2X}{\partial t\partial
x^q}\Big)+2\Big(\dfrac{\partial X}{\partial
x^q},\dfrac{\partial^2X}{\partial t\partial x^l}\Big)\Big]\\
&~~-2g^{kl}\Gamma^m_{ij}\Big(\vec n,\dfrac{\partial^2 X}{\partial
t\partial x^k}\Big)\Big(\dfrac{\partial X}{\partial
x^l},\dfrac{\partial^2 X}{\partial t\partial
x^m}\Big)-2\dfrac{\partial \Gamma^k_{ij}}{\partial t}\Big(\vec
n,\dfrac{\partial^2 X}{\partial t\partial
x^k}\Big)\\
&~~+2g^{kl}h_{ij}\Big(\vec n,\dfrac{\partial^2 X}{\partial t\partial
x^k}\Big)\Big(\vec n,\dfrac{\partial^2 X}{\partial t\partial
x^l}\Big)+\Big(\vec n,\dfrac{\partial}{\partial
x^i}\left(\dfrac{\partial H}{\partial x^j}\vec
n-Hh_{jk}g^{kl}\dfrac{\partial X}{\partial
x^l}\right)\Big)\\
&=\nabla_i\nabla_jH-Hh_{jk}g^{kl}h_{il}+g^{kl}h_{ij}\Big(\vec
n,\dfrac{\partial^2 X}{\partial t\partial x^k}\Big)\Big(\vec
n,\dfrac{\partial^2 X}{\partial t\partial
x^l}\Big)-2\dfrac{\partial\Gamma^k_{ij}}{\partial t}\Big(\vec
n,\dfrac{\partial^2 X}{\partial t\partial
x^k}\Big).\\\end{aligned}\end{equation*} Using (\ref{5.1}), we
obtain$$\dfrac{\partial^2h_{ij}}{\partial t^2}=\bigtriangleup
h_{ij}-2Hg^{kl}h_{il}h_{jk}+|A|^2h_{ij}+g^{kl}h_{ij}\Big(\vec
n,\dfrac{\partial^2X}{\partial t\partial x^k}\Big)\Big(\vec
n,\dfrac{\partial^2X}{\partial t\partial
x^l}\Big)-2\dfrac{\partial\Gamma^k_{ij}}{\partial t}\Big(\vec
n,\dfrac{\partial^2X}{\partial t\partial x^k}\Big).
$$ This proves (\ref{5.5}). Thus, the proof of Theorem 5.1 is
completed  .$\quad\quad\blacksquare$

\begin{Theorem} Under the hyperbolic mean curvature
flow,\begin{equation}\label{5.6}\begin{aligned}\dfrac{\partial^2H}{\partial
t^2} &=\bigtriangleup
H+H|A|^2-2g^{ik}g^{jl}h_{ij}\Big(\dfrac{\partial^2X}{\partial
t\partial x^k},\dfrac{\partial^2X}{\partial t\partial
x^l}\Big)+Hg^{kl}\Big(\vec n,\dfrac{\partial^2X}{\partial t\partial
x^k}\Big)\Big(\vec n,\dfrac{\partial^2X}{\partial t\partial
x^l}\Big)\\
&~~-2g^{ij}\dfrac{\partial\Gamma^k_{ij}}{\partial t}\Big(\vec
n,\dfrac{\partial^2X}{\partial t\partial
x^k}\Big)+2g^{ik}g^{jp}g^{lq}h_{ij}\dfrac{\partial g_{pq}}{\partial
t}\dfrac{\partial g_{kl}}{\partial t}-2g^{ik}g^{jl}\dfrac{\partial
g_{kl}}{\partial t}\dfrac{\partial h_{ij}}{\partial
t},\end{aligned}\end{equation}
\begin{equation}\label{5.7}\begin{aligned}\dfrac{\partial^2}{\partial t^2}|A|^2
&=\bigtriangleup(|A|^2)-2|\nabla A|^2+2|A|^4+2|A|^2g^{pq}\Big(\vec
n,\dfrac{\partial^2X}{\partial t\partial x^p}\Big)\Big(\vec
n,\dfrac{\partial^2X}{\partial t\partial
x^q}\Big)\\
&~~+2g^{ij}g^{kl}\dfrac{\partial h_{ik}}{\partial t}\dfrac{\partial
h_{jl}}{\partial t}-8g^{im}g^{jn}g^{kl}h_{jl}\dfrac{\partial
g_{mn}}{\partial t} \dfrac{\partial h_{ik}}{\partial
t}-4g^{im}g^{jn}g^{kl}h_{ik}h_{jl}\Big(\dfrac{\partial^2X}{\partial
t\partial x^m},\dfrac{\partial^2X}{\partial t\partial
x^n}\Big)\\
&~~+2g^{im}\dfrac{\partial g_{pq}}{\partial t} \dfrac{\partial
g_{mn}}{\partial
t}h_{ik}h_{jl}\Big(2g^{jp}g^{nq}g^{kl}+g^{jn}g^{kp}g^{lq}\Big)-4g^{ij}g^{kl}h_{jl}\dfrac{\partial\Gamma^p_{ik}}{\partial
t}\Big(\vec n,\dfrac{\partial^2X}{\partial t\partial
x^p}\Big).\end{aligned}\end{equation}
\end{Theorem} {\bf Proof.} Noting
$$g^{hm}g_{ml}=\delta^h_l,$$ we get
$$\dfrac{\partial g^{ij}}{\partial t}=-g^{ik}g^{jl}\dfrac{\partial g_{kl}}{\partial t},$$
$$\dfrac{\partial^2 g^{ij}}{\partial t^2}=2g^{ik}g^{jp}g^{lq}\dfrac{\partial g_{pq}}{\partial t}
\dfrac{\partial g_{kl}}{\partial t}-g^{ik}g^{jl}\dfrac{\partial^2
g_{kl}}{\partial t^2}.$$By a direct calculation, we have
\begin{equation*}\begin{aligned}\dfrac{\partial^2H}{\partial
t^2}&=\dfrac{\partial^2g^{ij}}{\partial t^2}h_{ij}+2\dfrac{\partial
g^{ij}}{\partial t}\dfrac{\partial h_{ij}}{\partial
t}+g^{ij}\dfrac{\partial^2h_{ij}}{\partial
t^2}\\
&=\Big(2g^{ik}g^{jp}g^{lq}\dfrac{\partial g_{pq}}{\partial t}
\dfrac{\partial g_{kl}}{\partial t}-g^{ik}g^{jl}\dfrac{\partial^2
g_{kl}}{\partial t^2}\Big)h_{ij}-2g^{ik}g^{jl}\dfrac{\partial
g_{kl}}{\partial t}\dfrac{\partial h_{ij}}{\partial
t}+g^{ij}\dfrac{\partial^2 h_{ij}}{\partial
t^2}\\
&=2g^{ik}g^{jp}g^{lq}h_{ij}\dfrac{\partial g_{pq}}{\partial
t}\dfrac{\partial g_{kl}}{\partial t}-2g^{ik}g^{jl}\dfrac{\partial
g_{kl}}{\partial t}\dfrac{\partial h_{ij}}{\partial
t}-g^{ik}g^{jl}h_{ij}\Big[-2Hh_{kl}+2\Big(\dfrac{\partial^2X}{\partial
t\partial x^k},\dfrac{\partial^2X}{\partial t\partial
x^l}\Big)\Big]\\
&~~+g^{ij}\Big[\nabla_i\nabla_jH-Hh_{il}h_{jk}g^{lk}+g^{kl}h_{ij}\Big(\vec
n,\dfrac{\partial^2X}{\partial t\partial x^k}\Big)\Big(\vec
n,\dfrac{\partial^2X}{\partial t\partial
x^l}\Big)-2\dfrac{\partial\Gamma^k_{ij}}{\partial t}\Big(\vec
n,\dfrac{\partial^2X}{\partial t\partial
x^k}\Big)\Big]\\
&=\bigtriangleup
H+H|A|^2-2g^{ik}g^{jl}h_{ij}\Big(\dfrac{\partial^2X}{\partial
t\partial x^k},\dfrac{\partial^2X}{\partial t\partial
x^l}\Big)+Hg^{kl}\Big(\vec n,\dfrac{\partial^2X}{\partial t\partial
x^k}\Big)\Big(\vec n,\dfrac{\partial^2X}{\partial t\partial
x^l}\Big)\\
&~~-2g^{ij}\dfrac{\partial\Gamma^k_{ij}}{\partial t}\Big(\vec
n,\dfrac{\partial^2X}{\partial t\partial
x^k}\Big)+2g^{ik}g^{jp}g^{lq}h_{ij}\dfrac{\partial g_{pq}}{\partial
t}\dfrac{\partial g_{kl}}{\partial t}-2g^{ik}g^{jl}\dfrac{\partial
g_{kl}}{\partial t}\dfrac{\partial h_{ij}}{\partial
t}.\end{aligned}\end{equation*}This is nothing but the desired
(\ref{5.6}).

On the other hand, by the definition of $|A|^2$ and the formula
(\ref{5.2}), a direct calculation gives
\begin{equation*}\begin{aligned}\qquad\dfrac{\partial^2}{\partial t^2}|A|^2
&=2\dfrac{\partial^2 g^{ij}}{\partial
t^2}g^{kl}h_{ik}h_{jl}+2\dfrac{\partial g^{ij}}{\partial
t}\dfrac{\partial g^{kl}}{\partial t}h_{ik}h_{jl}+8\dfrac{\partial
g^{ij}}{\partial t}g^{kl}\dfrac{\partial h_{ik}}{\partial
t}h_{jl}\\
&~~+2g^{ij}g^{kl}\dfrac{\partial^2h_{ik}}{\partial
t^2}h_{jl}+2g^{ij}g^{kl}\dfrac{\partial h_{ik}}{\partial
t}\dfrac{\partial h_{jl}}{\partial t}\\
&=2\Big(2g^{im}g^{jp}g^{nq}\dfrac{\partial g_{pq}}{\partial t}
\dfrac{\partial g_{mn}}{\partial t}-g^{im}g^{jn}\dfrac{\partial^2
g_{mn}}{\partial
t^2}\Big)g^{kl}h_{ik}h_{jl}\\
&~~+2g^{im}g^{jn}g^{kp}g^{lq}h_{ik}h_{jl}\dfrac{\partial
g_{mn}}{\partial t}\dfrac{\partial g_{pq}}{\partial
t}-8g^{im}g^{jn}\dfrac{\partial g_{mn}}{\partial
t}g^{kl}\dfrac{\partial h_{ik}}{\partial
t}h_{jl}+2g^{ij}g^{kl}\dfrac{\partial h_{ik}}{\partial
t}\dfrac{\partial h_{jl}}{\partial t}\\
&~~+2g^{ij}g^{kl}h_{jl}\Big[\nabla_i\nabla_kH-Hh_{ip}h_{kq}g^{pq}+g^{pq}h_{ik}\Big(\vec
n,\dfrac{\partial^2X}{\partial t\partial x^p}\Big)\Big(\vec
n,\dfrac{\partial^2X}{\partial t\partial
x^q}\Big)-2\dfrac{\partial\Gamma^p_{ik}}{\partial t}\Big(\vec
n,\dfrac{\partial^2X}{\partial t\partial
x^p}\Big)\Big]\\
&=4g^{im}g^{jp}g^{nq}g^{kl}\dfrac{\partial g_{pq}}{\partial t}
\dfrac{\partial g_{mn}}{\partial
t}h_{ik}h_{jl}-2g^{im}g^{jn}\Big[-2Hh_{mn}+2\Big(\dfrac{\partial^2X}{\partial
t\partial x^m},\dfrac{\partial^2X}{\partial t\partial
x^n}\Big)\Big]g^{kl}h_{ik}h_{jl}\\
&~~+2g^{im}g^{jn}g^{kp}g^{lq}h_{ik}h_{jl}\dfrac{\partial
g_{pq}}{\partial t} \dfrac{\partial g_{mn}}{\partial
t}-8g^{im}g^{jn}g^{kl}\dfrac{\partial g_{mn}}{\partial t}
\dfrac{\partial h_{ik}}{\partial
t}h_{jl}\\
&~~+2g^{ij}g^{kl}h_{jl}\nabla_i\nabla_kH-2H{\bf
tr}(A^3)+2g^{pq}|A|^2\Big(\vec n,\dfrac{\partial^2X}{\partial
t\partial x^p}\Big)\Big(\vec n,\dfrac{\partial^2X}{\partial
t\partial x^q}\Big)\\
&~~-4g^{ij}g^{kl}h_{jl}\dfrac{\partial\Gamma^p_{ik}}{\partial
t}\Big(\vec n,\dfrac{\partial^2X}{\partial t\partial
x^p}\Big)+2g^{ij}g^{kl}\dfrac{\partial h_{ik}}{\partial
t}\dfrac{\partial h_{jl}}{\partial t}\\
&=\bigtriangleup(|A|^2)-2|\nabla A|^2+2|A|^4+2|A|^2g^{pq}\Big(\vec
n,\dfrac{\partial^2X}{\partial t\partial x^p}\Big)\Big(\vec
n,\dfrac{\partial^2X}{\partial t\partial
x^q}\Big)\\
&~~+2g^{ij}g^{kl}\dfrac{\partial h_{ik}}{\partial t}\dfrac{\partial
h_{jl}}{\partial t}-8g^{im}g^{jn}g^{kl}h_{jl}\dfrac{\partial
g_{mn}}{\partial t} \dfrac{\partial h_{ik}}{\partial
t}-4g^{im}g^{jn}g^{kl}h_{ik}h_{jl}\Big(\dfrac{\partial^2X}{\partial
t\partial x^m},\dfrac{\partial^2X}{\partial t\partial
x^n}\Big)\\
&~~+2g^{im}\dfrac{\partial g_{pq}}{\partial t} \dfrac{\partial
g_{mn}}{\partial
t}h_{ik}h_{jl}\Big(2g^{jp}g^{nq}g^{kl}+g^{jn}g^{kp}g^{lq}\Big)-4g^{ij}g^{kl}h_{jl}\dfrac{\partial\Gamma^p_{ik}}{\partial
t}\Big(\vec n,\dfrac{\partial^2X}{\partial t\partial
x^p}\Big).\end{aligned}\end{equation*}This proves
(\ref{5.7}).$\qquad\qquad\qquad\qquad\blacksquare$

\section{Relations between hyperbolic mean curvature flow and the equations for
extremal surfaces in the Minkowski space $\mathbb{R}^{1,n}$}

In this section, we study the relations between the hyperbolic mean
curvature flow and the equations for extremal surfaces in the
Minkowski space $\mathbb{R}^{1,n}$.

Let $v=(v_0,v_1,\cdots,v_n)$ be a position vector of a point in the
$(1+n)$-dimensional Minkowski space $\mathbb{R}^{1,n}.$ The scalar
product of two vectors $v$ and $w=(w_0,w_1,\cdots,w_n)$ is$$\langle
v, w\rangle=-v_0w_0+\sum_{i=1}^{n}v_iw_i.$$ The Lorentz metric of
$\mathbb{R}^{1,n}$ reads $$ds^2=-dt^2+\sum_{i=1}^{n}(dx^i)^2.$$

A massless $n$-dimensional surface moving in $(1+n)$-dimensional
Minkowski space can be defined by letting its action be proportional
to the $(1+n)$-dimensional volume swept out in the Minkowski space.
It is a natural generalization of the massless string theory, and it
is interesting in its own right, as an example in which geometry,
classical relativity and quantum mechanics are deeply connected.
Hoppe et al \cite{B}, \cite{Ho} and Huang and Kong \cite{hk} have
obtained some interesting results about it.

We are interested in the following motion of an $n$-dimensional
Riemannian manifold in $\mathbb{R}^{1,n+1}$ with the following
parameter
\begin{equation}\label{6.1}(t,x^1,\cdots,x^n)\rightarrow \hat X=(t,~X(t,x^1,\cdots,x^n)),\end{equation}
where $(x^1,\cdots,x^n)\in \mathscr{M}$ and $\hat X(\cdot,t)$ be a
positive vector of a point in the Minkowski space
$\mathbb{R}^{1,n+1}$. The induced Lorentz metric reads
\begin{equation}\begin{cases}~\hat g_{00}=-1+\Big(\dfrac{\partial X}{\partial t},\dfrac{\partial X}
{\partial t}\Big),\vspace{2mm}\\~ \hat g_{0i}=\hat
g_{i0}=\Big(\dfrac{\partial X}{\partial t},\dfrac{\partial
X}{\partial x^i}\Big), \qquad (i, j=1,\cdots,n).\vspace{2mm}\\~\hat
g_{ij}=g_{ij}=\Big(\dfrac{\partial X}{\partial x^i},\dfrac{\partial
X}{\partial x^j}\Big)\end{cases}\end{equation}

By the variational method or by vanishing mean curvature of the
sub-manifold $\mathscr{M},$ we can obtain the following equation for
the motion of $\mathscr{M}$ in the Minkowski space
$\mathbb{R}^{1,n+1}$
\begin{equation}\label{6.2}\hat g^{\alpha\beta}\nabla_\alpha\nabla_\beta\hat X=\hat g^{\alpha\beta}\Big(
\dfrac{\partial^2\hat X} {\partial x^\alpha\partial
x^\beta}-\hat\Gamma_{\alpha\beta}^{\gamma}\dfrac{\partial \hat
X}{\partial x^\gamma}\Big)=0,\end{equation} where
$\alpha,\beta=0,1,\cdots,n.$ It is convenient to fix the
parametrization partially (see Bordemann and Hoppe \cite{B}) by
requiring
\begin{equation}\label{6.3}\hat g_{oi}=\hat g_{i0}=\Big(\dfrac{\partial
X}{\partial t},\dfrac{\partial X}{\partial
x^i}\Big)=0.\end{equation}It is easy to see that the equation
(\ref{6.2}) is equivalent to the following system
\begin{equation}\label{6.4}\Big(\dfrac{\partial^2X}{\partial t^2},\dfrac{\partial X}{\partial
t}\Big)-g^{ij}\Big(|X_t|^2-1 \Big)\Big(\dfrac{\partial^2X}{\partial
t\partial x^i},\dfrac{\partial X}{\partial
x^j}\Big)=0,\end{equation}
\begin{equation}\label{6.5}\begin{aligned}&\dfrac{\partial^2X}{\partial
t^2}+ g^{ij}\Big(\dfrac{\partial^2X}{\partial x^i\partial
x^j}-\Gamma_{ij}^{k}\dfrac{\partial X}{\partial
x^k}\Big)\Big(|X_t|^2-1\Big)-\dfrac{1}{|X_t|^2-1}
\left(\dfrac{\partial^2X}{\partial t^2},\dfrac{\partial X}{\partial
t}\right)\dfrac{\partial X}{\partial
t}\\&~~~~~~~~~~~~~+g^{kl}\left(\dfrac{\partial^2X}{\partial
t\partial x^l},\dfrac{\partial X}{\partial t}\right)\dfrac{\partial
X}{\partial x^k}+g^{ij}\left(\dfrac{\partial^2X}{\partial t\partial
x^i},\dfrac{\partial X}{\partial x^j}\right)\dfrac{\partial
X}{\partial t}=0,\end{aligned}\end{equation}where
$$|X_t|^2=\Big(\dfrac{\partial X}{\partial t},\dfrac{\partial
X}{\partial t}\Big).$$ We observe that,  when $\frac{\partial
X}{\partial t}\rightarrow 0,$ the limit of the equation (\ref{6.4})
reads
\begin{equation}\label{6.6}g^{ij}\Big(\dfrac{\partial^2X}{\partial t\partial x^i},\dfrac{\partial X}{\partial x^j}\Big)=0
,\end{equation}i.e.,
\begin{equation}\label{6.7}\dfrac{\partial}{\partial t}\det(g_{ij})=0.\end{equation} Moreover, the
equation (\ref{2.2}) is nothing but the limit of the equation
(\ref{6.5}) as $\frac{\partial X}{\partial t}$ approaches to zero.

\vskip 10mm\noindent{\Large {\bf Acknowledgements.}} The authors
thank the referee for his/her pertinent comments and valuable
suggestions. C.-L. He would like to thank the Center of Mathematical
Sciences at Zhejiang University for the great support and
hospitality. The work of Kong was supported in part by the NNSF of
China (Grant No. 10671124), the NCET of China (Grant No.
NCET-05-0390) and the Qiu-Shi Professor Fellowship from Zhejiang
University, China; the work of Liu was supported in part by the NSF
and NSF of China.

\end{document}